\documentclass[11pt,a4paper]{amsart}
\usepackage[utf8]{inputenc} 
\usepackage{url}
\usepackage{hyperref}
\usepackage{amsmath, amscd, amsfonts, amsthm, amssymb, latexsym, comment, stmaryrd}
\usepackage{faktor} 
\usepackage[british]{babel}
\usepackage{amsthm}
\usepackage{graphicx}
\usepackage{color}
\usepackage[top=2cm,bottom=2cm,left=2cm,right=2cm]{geometry}
\usepackage[all]{xy}
\usepackage[titletoc]{appendix}
\usepackage[T1]{fontenc}
\usepackage{fullpage}
\usepackage{verbatim}
\usepackage{chngcntr}
\newtheorem{thm}{Theorem}[section]
\newtheorem{lem}[thm]{Lemma}
\newtheorem{cor}[thm]{Corollary}

\newtheorem{pro}[thm]{Proposition}
\newtheorem{exa}[thm]{Example}
\newtheorem{defn}[thm]{Definition}
\newtheorem{rem}[thm]{Remark}
\newtheorem{ques}[thm]{Question}

\makeatletter \@addtoreset{thm}{subsection} \renewcommand{\thethm}{\ifnum\value{subsection}=0\relax\arabic{section}.\arabic{thm}\else\arabic{section}.\arabic{subsection}.\arabic{thm}\fi}  \makeatother

\newcommand{\A}{\alpha}

\newcommand{\gal}{\operatorname{Gal}}
\newcommand{\gl}{\operatorname{GL}}

\newcommand{\hm}{\operatorname{Hom}}

\newcommand{\im}{\operatorname{Im}}

\newcommand{\id}{\operatorname{id}}

\newcommand{\aut}{\operatorname{Aut}}
\newcommand{\en}{\operatorname{End}}

\newcommand{\function}[5]{%
\begin{array}{clcl}
  #1 : & #2  & \longrightarrow & #3\\
       & #4  & \longmapsto     & #5
  \end{array}
}

\title{Hopf-Galois Structures on  Extensions Arising from Division Points on Elliptic Curves}
\author{ Leolin Nkuete   }
\date{\today }

\begin{document}

\begin{abstract}
We provide a general construction of infinite Hopf-Galois extensions from a direct system of finite Hopf-Galois extensions.  We then apply this framework to extensions generated by division points on elliptic curves. Specifically, for a number  field  $K$ and an elliptic curve  $E/K$  of rank at least $1$ with $P\in E(K)$ of infinite order, we consider for each integer $n\geq 2$  the extension $L_n=K(Q_n)$, where $Q_n\in E(\overline{K})$ is such that $nQ_n=P$. We show that the finite-level extensions $L_n/K$ need not be  Hopf-Galois, and we give sufficient conditions for them to be so. 
We then  prove that the infinite extension $L_\infty=\bigcup_{n\geq 2}L_n$ of $K$ is potentially Hopf-Galois, meaning it becomes Hopf-Galois after a finite extension of the base field.
\end{abstract}

 \maketitle
\section{Introduction}

A finite field extension $L/K$ is a \emph{Hopf-Galois extension} if there exists a finite cocommutative  $K$-Hopf algebra $H$  such that $L$ is a $H$-module algebra and the $K$-linear map $$\operatorname{\underline{can}}:L\otimes_K H\longrightarrow \operatorname{End}_K(L),\, \alpha \otimes h\mapsto \operatorname{\underline{can}}(\alpha \otimes h)(\beta):= \alpha (h\cdot\beta)$$ is bijective. The pair $(H,\operatorname{\underline{can}})$ is called a Hopf-Galois structure on $L/K$. This notion is a natural generalisation of classical Galois theory, in the sense that every finite Galois extension is Hopf-Galois. Indeed, if $L/K$ is Galois with Galois group $\Gamma$, then one may take $H$ to be the group algebra $K[\Gamma]$. However, not all Hopf-Galois extensions are Galois (see \cite{Greither}). For separable extensions, Greither and Pareigis \cite{Greither} give the following group-theoretic characterisation of the existence of Hopf-Galois structures.
\begin{thm}\label{thm:G-P}
 Let $L/K$ be a separable finite field extension of degree $d$. Let $\widetilde{L}$ be its Galois closure, $G=\gal(\widetilde{L}/K)$ and $G'=\gal(\widetilde{L}/L)$. There  exists a bijective correspondence between the set of Hopf-Galois structures on $L/K$ and the set of regular subgroups of the symmetric group $S_d$ normalised by $\lambda(G)$, where $\lambda:G\rightarrow \operatorname{Perm}(G/G')$ is the morphism given by the action of $G$ on the set of left cosets $G/G'$ and $S_d$ is identified with $\operatorname{Perm}(G/G')$.
\end{thm}
Note that the homomorphism $\lambda$ in Theorem~\ref{thm:G-P} is injective, and hence, we may identify $G$ with its image under $\lambda$. Moreover, if the corresponding regular subgroup $N$ is contained in $G$, equivalently, if $G'$ admits a normal complement in $G$, then the extension $L/K$ is said to be \emph{almost classically Galois}. This is an important class of Hopf-Galois extensions, since in this case there is a bijective correspondence between the set of $K$-sub-Hopf algebras of $H$ and the set of intermediate fields of $L/K$. This gives a natural generalisation of the classical Galois correspondence between intermediate extensions and subgroups of the Galois group.
There is an equivalent reformulation of Theorem~\ref{thm:G-P} due to Childs \cite[Proposition 1]{Childs}, involving the holomorph $\operatorname{Hol}(N)$ of the regular subgroup $N$.

\begin{thm}\label{thm: Hol}
The Galois group $G$ normalises a regular subgroup $N$ of $\operatorname{Perm}(G/G')$ if and only if $G$ is a subgroup of the holomorph $\operatorname{Hol}(N):=N\rtimes\aut(N)$ of $N$.
\end{thm}
Since the $\operatorname{Hol}(N)\subset \operatorname{Perm}(G/G')$, in practice, it is more efficient  to check that $G$ is a subgroup of $\operatorname{Hol}(N)$ rather than to find a regular subgroup normalised by $G$.

The theory of finite Hopf-Galois extensions was extended to infinite extensions in \cite{Bui}, and the definition is generalised as follows.

\begin{defn}\label{def: infHopf}
Let $L/K$ be a field extension and $H$ be a proartinian  $K$-Hopf algebra such that $L$ is a discrete left $H$-module algebra. We say that the extension $L/K$ is Hopf-Galois or $H$-Galois  if the following canonical map defined in Section~\ref{subsec: infi},  $\operatorname{\textbf{\underline{can}}}: L\hat{\otimes}H \rightarrow \operatorname{End}_K(L)$ is an homeomorphism of topological vector spaces, where $L\hat{\otimes}H$ and $\operatorname{End}_K(L)$ are proartinian topological vector spaces.

\end{defn}

However, this theory lacks of examples, that is, there are very few known examples of infinite Hopf-Galois extensions that are not classical infinite Galois extensions. The aim of this paper is to provide a general method for constructing infinite Hopf-Galois extensions from an infinite family of finite Hopf-Galois extensions and provide explicit examples.  For this purpose, we consider $K$ to be a field and   fix an algebraic closure $\overline{K}$ of $K$. Let $(I,\leq)$ be a directed  partially ordered set and  let $(L_n)_{n\in I}$ be a family of finite separable extensions of $K$ such that $L_n\subseteq L_m$ whenever $n\leq m$. For every $n\in I$, we denote by  $\widetilde{L}_n$  the Galois closure of $L_n/K$ in $\overline{K}$ and we  let $$G_n:=\gal(\widetilde{L}_n/K),\qquad G_n'=\gal(\widetilde{L}_n/L_n)\quad \text{and}\quad X_n=\hm_{K\operatorname{-Alg}}(L_n,\widetilde{L}_n)\simeq G_n/G_n'.$$ If $n\leq m$, then $\widetilde{L}_n\subseteq\widetilde{L}_m$, and  the restriction  homomorphism $r_{n,m}:G_m\longrightarrow G_n$ is surjective.  Let $$\psi_{n,m}: X_m\rightarrow X_n,\; x\mapsto x_{|L_n}$$ be the restriction map. Suppose that for all $n\in I$, the extension $L_n/K$ is Hopf-Galois for a finite-dimensional $K$-Hopf algebra $H_n$ and let $N_n\subset \operatorname{Perm}(X_n)$ be its Greither-Pareigis group, where the action of $N_n$ on $X_n$ is on the right, that is, $x\cdot \eta:=\eta^{-1}(x)$ for all $\eta\in N_n$ and $x\in X_n$. The $K$-Hopf algebra $H_n$ is an $\widetilde{L}_n$-form of $K[N_n]$ obtained via Galois descent. As the $\widetilde{L}_n$-forms of $K[N_n]$ are classified by $ \operatorname{H}^1(G_n,\aut(N_n))=\hm_{\operatorname{Group}}(G,\aut(N))$ (see Theorem~\ref{them: Galois_class} and Remark~\ref{rem: forms}), we denote by $\theta^n\in \operatorname{H}^1(G_n,\aut(N_n))$ the associated cocycle. Note that for any $\sigma\in G_n$ and $\eta\in N_n$, we have $\theta_\sigma^n(\eta)=\sigma\eta\sigma^{-1}$ viewed in $\operatorname{Perm}(X_n)$ (see \cite[page 248]{Greither}). In particular, when $L_n/K$ is almost classically Galois, the conjugation is viewed in $G_n$.

\begin{thm}\label{thm: main}
With the notation as above, let $L=\bigcup_{n\in I}L_n.$ Suppose that  for every $n\leq m$, there exists a surjective  group homomorphism $ q_{n,m}:N_m\longrightarrow N_n $ satisfying the following conditions:
\begin{enumerate} 
\item \label{cond:1} $q_{n,n}=\id_{N_n}$ and $q_{l,m}=q_{l,n}\circ q_{n,m}$ whenever $l\leq n\leq m$. 
\item \label{cond: 2} For every $\sigma\in G_m$, we have $ q_{n,m}\circ\theta^m_\sigma = \theta^n_{r_{n,m}(\sigma)}\circ q_{n,m}$ whenever $n\leq m$. 
\item \label{cond: 3} For every $\eta\in N_m$ and every $x\in X_m$, $ \psi_{n,m}(x\cdot\eta) =  \psi_{n,m}(x)\cdot q_{n,m}(\eta).$
\end{enumerate} Then $q_{n,m}$ induces a surjective  morphism of $K$-Hopf algebras $\pi_{n,m}:H_m\longrightarrow H_n$.  These morphisms satisfy $\pi_{n,n}=\id_{H_n}$ and $\pi_{l,m}=\pi_{l,n}\circ\pi_{n,m}$ whenever $l\leq n\leq m$. Consequently, $ H:=\varprojlim_n H_n$ is a proartinian $K$-Hopf algebra and the extension $L/K$ is $H$-Galois.
\end{thm}

Condition ~(\ref{cond: 2}) ensures that the morphism $\widetilde{L}_m[N_m]\rightarrow \widetilde{L}_m[N_n]$ of $\widetilde{L}_m$-Hopf algebras  induced by $q_{n,m}$ descends to a morphism $\pi_{n,m}$ of the corresponding $K$-Hopf algebras. Condition~(\ref{cond: 3}) ensures that the Hopf-Galois actions at the finite levels are compatible, that is, for any $y\in L_n$ and $h_m\in H_m$, we have $h_m\cdot y=\pi_{n,m}(h_m)\cdot y$ (see Lemma \ref{lem: equivaction}). 

To construct examples satisfying Theorem~\ref{thm: main}, we assume that $K$ is a number field. Let $E/K$ be an elliptic curve of rank at least $1$, and let $P\in E(K)$ be a fixed rational point of infinite order. For any integers  $n,m\geq 2$, let $Q_n, Q_m\in E(\overline{K})$ be such that $nQ_n = P$, and $\frac{m}{n}Q_m=Q_n$ whenever $n\mid m$. For any integer $n\geq 2$, we consider the field $L_n:=K(Q_n),$ obtained by adjoining to $K$ the coordinates of $Q_n$. Clearly, we have $L_n\subset L_m$ whenever $n\mid m$.  One might first ask the following question.

\begin{ques}\label{question}
Do the extensions $L_n/K$ admit Hopf-Galois structures for all $n\geq 2$? 
\end{ques}
Although the main motivation is to construct an infinite Hopf-Galois extension, Question~\ref{question} is also of independent interest. We therefore investigate the existence of Hopf-Galois structures on the extensions $L_n/K$, which are generally not Galois because $L_n$ need not contain the $n$-torsion points of $E$. This is analogous to the radical extensions $K(\sqrt[n]{a})/K$, which need not be Galois when $K$ does not contain the $n$-th roots of unity. 

Controlling the degree $[L_n:K]$ gives essential information on whether the extensions $L_n/K$ can be Hopf-Galois. These extensions are ``Kummer-like'' and therefore also belong to a broader framework generalizing classical Kummer theory to settings in which the base field need not contain the relevant $n$-torsion points or the relevant roots of unity in the radical case (see \cite{perucca}).

To study Question~\ref{question}, for each integer $n\geq 2$, we let $G_n:=\gal(\widetilde{L}_n/K)$. We consider  the map $$\delta_n: G_n\longrightarrow E[n]\,, \sigma\mapsto  \sigma(Q_n)-Q_n.$$ This map is a $1$-cocycle of $G_n$ with values in $E[n]$. Considering the $G_n$-stable set  $$I_n=\{\sigma\delta_n(\tau)\;|\, \sigma,\tau\in G_n\}= \{\sigma\tau(Q_n)-\sigma(Q_n)\;|\, \sigma,\tau\in G_n\},$$  the field $K_n:=K(I_n)$ is a Galois extension of $K$ such that the composite field $K_nL_n$ is equal to $\widetilde{L}_n$. We then prove in Proposition~\ref{prop: Im} that if the image of $\delta_n$ is a subgroup of $E[n]$, then the extension $L_n/K$ is Hopf-Galois. As a consequence, if   $K_n\cap L_n=K$, then $L_n/K$ is almost classically  Galois, where the corresponding Greither-Pareigis group is isomorphic to $\gal(\widetilde{L}_n/K_n)\cong\im\delta_n$ (see Corollary~\ref{cor: m=1}). These results lead to the following theorem.

\begin{thm}\label{Thm: main}

Let $K$ be a number field and let  $E/K$ be an elliptic curve of rank $1$. Let $P \in E(K)$ be a point of infinite order. For an integer $n\geq 2$, let $Q_n\in E(\overline{K})$ be such that $nQ_n = P$. Let $F_n:=K_n\cap L_n$. If $[F_n: K]\leq 2$, then the extension $L_n/K$ is Hopf-Galois.
\end{thm}

Theorem~\ref{Thm: main} suggests that counterexamples may arise when the degree $[F_n:K]$ is sufficiently large. However, a large intersection degree alone does not imply that $L_n/K$ is not Hopf-Galois. To apply the Greither-Pareigis correspondence and its holomorph reformulation, one must also control the degree $[L_n:K]$. Motivated by this observation, in Section~\ref{subsec: counter}, we give such a construction and  obtain the following result, which gives, under suitable conditions, information about the degree of the extension $L_n/K$.

\begin{thm}\label{thm: deg}

Let $K$ be a number field and let $E/K$ be an elliptic curve of rank at least $1$. Suppose that the mod-$\ell$ Galois representation associated with $E$ is surjective for all primes $\ell$ outside a finite set. Then there exist primes $p$ and $q$, and a point $P\in E(K)$ of infinite order, such that, for $n=pq$, there exists a point $Q_n\in E(\overline{K})$ satisfying $nQ_n=P$ and such that   $$[K(Q_n):K]\in\{(p^2-1)mq,\; 1\leq m\le q \}.$$  Moreover, $p$ divides the degree of the Galois closure $\widetilde{L}_n$ of $K(Q_n)$ over $K$.
\end{thm}

Using Theorem~\ref{thm: deg} and Theorem~\ref{thm: Hol}, we give a negative answer to Question \ref{question} by constructing an explicit example of an extension $L_n/K$  such that there is a prime number $p$ dividing $\#\gal(\widetilde{L}_n/K)$ but not dividing $\#\operatorname{Hol}(N)$ for any group $N$ of order $[L_n:K]$ (see Corollary~\ref{cor: counter}). 

Thus, in our setting, the obstruction to the existence of a Hopf-Galois structure lies in the intersection $F_n=K_n\cap L_n$. However, as we will see in Section~\ref{section 4}, after a finite extension of the base field, this obstruction uniformly disappears as it follows. For  $n\mid m$, we have  $I_n =\frac{m}{n} I_m$. This implies that  $K_n\subset K_m$, and therefore, $F_n\subset F_m$.  We let $K'=\bigcup_{n\geq 2}F_n$. As there exists a constant $C=C(E,P,K)$ such that $[F_n:K]\leq C$ for all $n$ (see \cite[Theorem 1]{Bertrand}), it follows that $K'$ is a finite extension of $K$ such that $K'K_n\cap K'L_n=K'$ for all $n\geq 2$.

\begin{defn}
A field extension (not necessarily finite) is called  potentially Hopf-Galois if it becomes Hopf-Galois after a finite extension of the base field.
\end{defn}

\begin{defn}
 An infinite Hopf-Galois extension is called almost classically Galois if it is a direct limit of finite almost classically Galois extensions.
\end{defn}

\begin{thm}\label{thm: main2}
Let $K$ be a number field and let $E/K$ be an elliptic curve of rank at least $1$ defined over $K$. Let $P \in E(K)$ be a point of infinite order. For any integers  $n,m\geq 2$, let $Q_n, Q_m\in E(\overline{K})$ such that $nQ_n = P$, and $\frac{m}{n}Q_m=Q_n$ whenever $n\mid m$.  Let $L_n := K(Q_n)$ and $L_\infty:=\bigcup_n L_n$. Then the extension $L_\infty/K$ is potentially almost classically Galois.
\end{thm}

This paper is organised as follows. In Section~\ref{section 2}, we recall the notion of an infinite Hopf-Galois extension and use the classification of $K$-forms of algebras to prove Theorem~\ref{thm: main}. In Section~\ref{section 3}, show that the extensions $K(Q_n)/K$  are not always Hopf-Galois and give  sufficient conditions for then to be so.  Finally, in Section~\ref{section 4}, we prove Theorem~\ref{thm: main2}. 

 \textbf{Notation:} If $S$ is a set, we denote $\#S$ for the cardinality of $S$. The following diagram gives a summary of other pieces of notation.
\[
\xymatrix@C=1.5em@R=1.5em{
& & \widetilde{L}_m  \ar@{-}[dll]_{\Gamma_m} \ar@{-}[dr]^{G_m'} \ar@{-}[dd] & \\
K_m \ar@{-}[ddr]\ar@{-}[rd]&  &    & \ar@{-}[dd] \ar@{-}[dll] L_m\\
& F_m \ar@{-}[ddr]   & \widetilde{L}_n=K_nL_n \ar@{-}[dd]^{\Gamma_n G_n'}\ar@{-}[dr]^{G_n'}  \ar@{-}[dl]_{\Gamma_n} &  \\
& K_n:=K(I_n) \ar@{-}[dr] &   & L_n:=K(Q_n) \ar@{-}[dl] \\
& & F_n:=K_n\cap L_n \ar@{-}[d] & \\
& & K &
}
\]

\section{A construction of infinite Hopf Galois extensions}\label{section 2}
In this section, we briefly introduce the theory of infinite Hopf-Galois extension (for more details see \cite{Bui}), and we prove Theorem~\ref{thm: main}.

\subsection{Infinite Hopf-Galois extensions}\label{subsec: infi}
 
\begin{defn}
A proartinian $K$-vector space $V$ is a topological $K$-vector space such that there exists a family $\Lambda_V$ of open subspaces $U\subseteq V$ that form a basis of open neighbourhoods of $0$ such that the following hold:
\begin{enumerate}
    \item the quotient $V/U$ is a finite dimensional discrete $K$-vector space for each $U\in\Lambda_V$,
    \item the natural $K$-linear map $\pi^V: V\longrightarrow  \varprojlim_{U\in \Lambda_V} V/U$ induced by all projection maps $\pi_U^V: V\rightarrow V/U$, is bijective.
\end{enumerate}
A morphism between proartinian $K$-vector spaces is a continuous $K$-linear map.
\end{defn}
We define the completed tensor product $\widehat{\otimes}$ between two proartinian $K$-vector spaces $V$ and $W$ by $V\widehat{\otimes} W:= \varprojlim_{(U_V,U_W)\in \Lambda_V\times \Lambda_W}V/U_V\otimes W/U_W$.

\begin{defn}
Let $H$ be a  proartinian $K$-vector space. We say that $H$ is a proartinian $K$-Hopf algebra if there is a basis $\Lambda_H$ of $0$ consisting of open subspaces $U$ of $H$ such that for every $U\in \Lambda_H$, the quotient $H/U$ is a finite dimensional $K$-Hopf algebra and for all $U,V\in\Lambda_H$ with $U\subseteq V$, the natural maps $\pi_{V,U}:H/U\longrightarrow H/V$ are morphisms of $K$-Hopf algebras with respect to the fixed Hopf algebra structures.
\end{defn}
Note that extra structure morphisms (multiplication, a unit, a counit and an antipode morphism) on $H$ can be constructed by taking a projective limit of the structure maps of all $H/U$ for $U\in\Lambda_H$.

\begin{exa}
Any finite dimensional  $K$-Hopf algebra is a discrete proartinian Hopf algebra.
\end{exa}

\begin{defn}
 Let $H$ be a proartinian $K$-Hopf algebra and $\Lambda_H$ a basis of neighbourhoods of $0$ for $H$. We say that a $K$-subspace $I\subseteq H$ is a proartinian Hopf ideal if $(I+U)/U\subseteq H/U$ is a Hopf ideal for all $U\in\Lambda_H$. 
\end{defn}

\begin{defn}
Let $H$ be a proartinian $K$-Hopf algebra. A proartinian $K$-vector space $L$ is called a proartinian  left $H$-module algebra if $L$ is an $H$-module with a continuous structure map $H\times L\rightarrow L$ and carries a compatible structure of proartinian algebra.
\end{defn}

Let $H$ be a proartinian $K$-Hopf algebra and  $L$ be a discrete $H$-module algebra. Then $\en_K(L)=\varprojlim_{I}\hm_K(L^I,L)$ is a proartinian $L$-vector space (see \cite[Lemma 3.18]{Bui}), where $I$ varies over all open proartinian Hopf ideals in $H$ and $L^I=\{x\in L\,|\,\forall i\in I, i\cdot x=0  \}$. The canonical continuous  map $\textbf{\underline{can}}: L\hat{\otimes}H \rightarrow \operatorname{End}_K(L)$ in  Definition~\ref{def: infHopf} is the unique morphism that makes the following diagram commutative for any open proartinian Hopf ideal $I\subset H$.  

$$\xymatrixcolsep{3.5pc}\xymatrix{
	L\widehat{\otimes} H  \ar[r]^-{\textbf{\underline{can}}} \ar[d]_{\pi_I} & \en_K(L) \ar[d]^{res_I}\\
	L\otimes H/I \ar[r]_-{\operatorname{\underline{can}}_I } &\hm_K(L^I,L)}$$
where $\operatorname{\underline{can}}_I(x\otimes\bar{h}):=x(\bar{h}\cdot y) $ for all $x\in L$, $y\in L^I$ and $\bar{h}\in H/I$.
 
Let $L/K$ be an $H$-Galois extension. Let $K\subset L_0\subseteq L$ be an intermediate field. We define the annihilator of $L_0$ in $H$ by $$J_H(L_0)=\{h\in H\,|\;\forall x\in L_0,\, h\cdot x=0 \}.$$
This is a  proartinian ideal of $H$ (see \cite[Lemma 3.19]{Bui}).
\begin{defn}
 Let $L/K$ be an $H$-Galois extension and $L_0$ be an intermediate field. We say that $L_0/K$ is an $H$-sub-Galois extension if $L_0$ is $H/J_H(L_0)$-Galois in the sense of Definition~\ref{def: infHopf}.  If Moreover, $H\cdot L_0\subset L_0$, we say that   $L_0$ is an $H$-normal sub-extension of $H$.
\end{defn}

\begin{thm}\cite[Theorem 2.24]{Bui}\label{thm: corespo}
 Let $L/K$ be an $H$-Galois extension and $L_0$ be an $H$-normal sub-extension. Then $L^{J_H(L_0)}=L_0$.
\end{thm}
 As a simple example of an infinite Hopf-Galois extension, if $L/K$ is a Galois extension (not necessarily finite) with Galois group $G$, then $L/K$ is Hopf-Galois for the completed group algebra $H = K\llbracket G\rrbracket$. The first example of an infinite Hopf-Galois extension that is not a classical infinite Galois extension is given in \cite[Example 3.5]{Bui}, where they consider the family of almost classical Galois extensions $(L_n:=\mathbb{Q}(\sqrt[3^n]{2}))_{n\in\mathbb{N}}$ and use the axiom of choice to obtain, for each $n\in\mathbb{N}$, a Hopf algebra $H_n$ turning $L_n$ into an $H_n$-Galois extension and such that for each $n,m\in\mathbb{N}$ with $m\leq n$, there is a surjective morphism $\pi_{m,n}:H_n\longrightarrow H_m$ leading to the Hopf-Galois extension $L=\bigcup_{n\in\mathbb{N}}L_n$ for the proartinian Hopf algebra $H=\varprojlim_n H_n$. In the rest of this section, we will use the classification of forms to generalise this construction to any suitable family of finite separable  Hopf-Galois extensions.

\subsection{Classification of $K$-forms}\label{subsec: 2.2}
In this subsection, we review some results on the classification of forms of algebras. For more details see \cite{Knustheorydeladescente}. Let $k$ be a ring and $K$  a commutative $k$-algebra.  We denote $\otimes$ the tensor product over $k$.

\begin{defn}
 Let $A$ be a $k$-algebra with possible extra algebraic structure (Hopf algebra). We say that a $k$-algebra $B$ with the same algebraic structure as $A$ is a $K$-form of $A$ if there exists an isomorphism of  $K$-algebras  $\beta: K\otimes B\longrightarrow K\otimes A$ which respects the extra algebraic structure.

 Two $K$-forms of $\beta: K\otimes B\longrightarrow K\otimes A$ and $\beta': K\otimes B'\longrightarrow K\otimes A$ are isomorphic if there exists an isomorphism of $k$-algebras $\phi:B\longrightarrow B'$ such that the following diagram commutes. $$ \xymatrixcolsep{3.5pc}\xymatrix{
            K\otimes B \ar[dr]_{\beta} \ar[r]^{1\otimes\phi}&  K\otimes B'\ar[d]^{\beta'} \\             
               &  K\otimes A  }$$
\end{defn}

Let us fix a $k$-algebra $A$. The next result tells us that when $K$ is a Galois extension of $k$ with Galois group $G$, then the $K$-forms of $A$ can be classified in terms of the $1$-cocycles of $G$ with values in $\aut_K(K\otimes A)$, the group of $K$-linear automorphisms that preserve the structure of $K\otimes A$. Here, $G$ acts on $K\otimes A$ via the left factor  and so,  $G$ acts on $\aut_{K}(K\otimes A)$ via $\sigma(\A)=\sigma\A\sigma^{-1}$ for any $\sigma\in G$ and $\A\in \aut_{K}(K\otimes A)$. If $\beta: K\otimes B\longrightarrow K\otimes A$ is a $K$-form of $A$ then the quantity $$ \theta_\sigma:=\theta(\sigma) = \beta\sigma\beta^{-1}\sigma^{-1},$$ which measures the failure of $\beta$ of being $G$-invariant  defines the $1$-cocycle associated to the isomorphism class of $(B,\beta)$ (see \cite[page 65]{Knustheorydeladescente}). Now, given $\theta\in \operatorname{H}^1(G,\, \aut_{K }(K\otimes A))$, we define for each $\sigma\in G$ the automorphism of $k$-algebras $\Bar{\sigma}:=\theta_\sigma\sigma: K\otimes A\longrightarrow K\otimes A$. It is simple to check that $\Bar{\sigma}$ is $K$-semi-linear and for any $\sigma,\tau\in G$, we have $$ \overline{\sigma\tau}=\theta_{\sigma\tau}\sigma\tau=\theta_\sigma\sigma\theta_\tau\sigma^{-1}\sigma\tau=\theta_\sigma\sigma\theta_\tau\tau=\Bar{\sigma}\Bar{\tau}.$$ Therefore, the family $\{\Bar{\sigma}\}_{\sigma\in G}$ is a Galois descent datum for $K\otimes A$ over $K$ (see \cite[Theorem 5.1]{Knustheorydeladescente}). We denote $B:=(K\otimes A)^G=\{x\in K\otimes A: \Bar{\sigma}(x)= x\}$ the descended $k$-algebra and  $\beta: K\otimes B\longrightarrow K\otimes A$ the corresponding isomorphism.  It follows that $(B,\beta)$ is a $K$-form of $A$ and, since by  definition $\beta$ satisfies $\Bar{\sigma}\beta=\beta\sigma$ for any $\sigma\in G$,  it follows that $\theta_\sigma=\beta\sigma\beta^{-1}\sigma^{-1}$. 

\begin{thm}\cite[Theorem 9.1]{Knustheorydeladescente}\label{them: Galois_class}
Suppose that  the extension $K/k$ is Galois with Galois group $G$. Then  the group 
 $  \operatorname{H}^1(G,\, \aut_{K }(K\otimes A))$  classifies the $K$-forms of $A$. 
\end{thm}
 This Theorem tells us that two $K$-forms of $A$ are isomorphic if and only if their associated cocycles are cohomologous.
\begin{defn}
Let $B_i$ be a $K$-form of $A_i$, $\beta_i: K\otimes B_i\longrightarrow K\otimes A_i$ ($i=1,2$). Then a morphism $f: K\otimes A_1\longrightarrow K\otimes A_2$ is called descendable if there is a morphism $g$ such that the diagram 
$$\xymatrix{ K\otimes B_1 \ar[r]^{1\otimes g} \ar[d]_{ \beta_1} & K\otimes B_2 \ar[d]^{\beta_2 } \\ K\otimes A_1\ar[r]_-{ f} & K\otimes A_2 }$$
commutes. If $g$ exists, it is unique.
\end{defn}

\begin{pro}\cite[Lemma 1.5]{Greither}{\label{pro: desct}}
 With the notation of the above definition, let $\theta^i$ be the associated cocycles of $B_i$. Then $f$ is descendable if and only if $$ f\theta^1_{\sigma}=\theta^2_{\sigma}f \quad\quad \text{ for all } \sigma\in G.$$
\end{pro}

\begin{rem}\cite[Theorem 3.1]{Greither}\label{rem: forms}
Let $L/K$ be a finite separable  $H$-Galois extension with the corresponding Greither-Pareigis group  $N$, and $\widetilde{L}$ be the Galois closure of $L/K$. Let $G=\gal(\widetilde{L}/K)$. Since the automorphisms of the $\widetilde{L}$-Hopf algebra  $\widetilde{L}\otimes K[N]=\widetilde{L}[N]$ are given by  automorphisms of $N$, it follows that the action of  $G$ on  $\aut_{\widetilde{L}\text{-Hopf}}(\widetilde{L}[N])$ by conjugation is trivial. Therefore, $\operatorname{H}^1(G,\aut_{\widetilde{L}\text{-Hopf}}(\widetilde{L}[N]))=\hm_{\operatorname{Group}}(G,\aut(N))$.  Using the regular subgroup $N$, Greither and Pareigis define the $1$-cocycle $\theta: G\rightarrow \aut(N)$  by $\theta_\sigma(\eta)=\sigma\eta\sigma^{-1}$ for all $\sigma\in G$, $\eta\in N$.  They obtain the Hopf algebra $H$ as the descended $K$-algebra. Thus, $H$ is an $\widetilde{L}$-form of $K[N]$.
\end{rem}

\subsection{  From finite Hopf-Galois extensions to infinite Hopf-Galois extension }

Let $(L_n)_{n\in I}$ be a family of finite separable Hopf-Galois extensions as given in the introduction.  For each $n\in I$, the $K$-Hopf algebra $H_n$ is an $\widetilde{L}_n$-form of $K[N_n]$ given by an isomorphism $\beta_n: \widetilde{L}_n\otimes H_n\rightarrow \widetilde{L}_n[N_n]$. Let $\theta^n\in \operatorname{H}^1(G_n,\, \aut_{\widetilde{L}_n }(\widetilde{L}_n[N_n]))=\operatorname{H}^1(G_n,\aut(N_n))$ be the associated cocycle. Let $n,m\in I$ such that $n\leq m$. The $K$-Hopf algebra $H_n$ becomes an $\widetilde{L}_m$-form of $K[N_n]$ via scalar extension $$  \widetilde{\beta}_n:=(1\otimes \beta_n): \widetilde{L}_m\otimes H_n\rightarrow \widetilde{L}_m[N_n],\; (a\otimes h)\mapsto a\beta_n(1\otimes h)$$ and the associated cocycle $\widetilde{\theta}^n\in \operatorname{H}^1(G_m,\, \aut_{\widetilde{L}_m }(\widetilde{L}_n[N_n])) $ is given by  $\widetilde{\theta}^n_\sigma:= \widetilde{\beta}_n \sigma \widetilde{\beta}_n^{-1}\sigma^{-1}$ for all $\sigma\in G_m$. The restriction homomorphism $r_{n,m}:G_m\to G_n$ induces an inflation  homomorphism $$\operatorname{Inf}:\operatorname{H}^1(G_n,\aut(N_n))\rightarrow \operatorname{H}^1(G_m,\aut(N_n))$$ defined by    $\operatorname{Inf}(\theta^n)(\sigma)=\theta^n_{r_{n,m}(\sigma)}$ for all $\sigma\in G_m$.  Since $\beta_n$ has coefficients in $\widetilde{L}_n$ and  the action of $G_m$ on   $\widetilde{L}_n$ factors through $r_{n,m}$, we have $$\widetilde{\theta}^n_\sigma= 1\otimes (\beta_n r_{n,m}(\sigma)\beta_n^{-1}r_{n,m}(\sigma)^{-1})= 1\otimes \theta^n_{r_{n,m}(\sigma)}.$$  After the  identifications $\aut_{\widetilde{L}_n }(\widetilde{L}_n[N_n])=\aut(N_n)=\aut_{\widetilde{L}_m }(\widetilde{L}_m[N_n])$, it follows that, at the level automorphism of $N_n$, we have $\widetilde{\theta}^n_\sigma= \theta^n_{r_{n,m}(\sigma)}$. Thus, $\operatorname{Inf}(\theta^n)$ and $\widetilde{\theta}^n$ correspond to the same $\widetilde{L}_m$-form $H_n$.


Now, for $i\in\{n,m\}$, we give the relation between the action of $N_i$ on $X_i$ and the natural Hopf-Galois action $\mu_i: H_i\otimes L_i\rightarrow L_i$ of $H_i$ on $L_i$. 
For any $i\in\{n,m\}$, there is an isomorphism of $\widetilde{L}_m$-algebras $$\phi_i: \widetilde{L}_m\otimes L_i\rightarrow \operatorname{Maps}(X_i,\widetilde{L}_m)=\widetilde{L}_m^{X_i}$$ defined by $\phi_i( a\otimes y)(x)=ax(y)$ for all $x\in X_i$ (see \cite[Lemma 2.1]{Greither}). The right action of $N_i$ on $X_i$ induces a left action of $N_i$ on $\widetilde{L}_m^{X_i}$ given by $(\eta\cdot\varphi)(x)=\varphi(x\cdot\eta)$ for all $\eta\in N_i$, $\varphi\in \widetilde{L}_m^{X_i}$ and $x\in X_i$. This action extends by linearity to an action of $\widetilde{L}_m[N_i]$ on $\widetilde{L}_m^{X_i}$. Then, under the isomorphism $\phi_i$,  $\widetilde{L}_m[N_i]$ acts linearly  on $\widetilde{L}_m\otimes L_i$.  The action of $N_i$ on $X_i$ and the action of $N_i$ on $\widetilde{L}_m\otimes L_i$ are related by the following. For all $y\in L_i$, $x\in X_i$ and $\eta\in N_i$, we have
$$\phi_i(\eta\cdot(1\otimes y))(x)=(\eta\cdot\phi_i(1\otimes y))(x)=\phi_i(1\otimes y)(x\cdot\eta)=(x\cdot\eta)(y).$$
 After a scalar extension of the natural Hopf-Galois action $\mu_i: H_i\otimes L_i\rightarrow L_i$, the following diagram commutes (see \cite[page 249]{Greither}).
$$ \xymatrixcolsep{2.5pc}\xymatrix{
            (\widetilde{L}_m\otimes H_i)\otimes_{\widetilde{L}_m} (\widetilde{L}_m\otimes L_i) \ar[d]_{1\otimes\mu_i
} \ar[r]^-{\beta_i\otimes 1 }&  \widetilde{L}_m[N_i]\otimes_{\widetilde{L}_m}  (\widetilde{L}_m\otimes L_i)\ar[d]^{}  \\  \widetilde{L}_m\otimes L_i  \ar@{=}[r]^{}     &  \widetilde{L}_m\otimes L_i  }   $$

This shows that, after extending scalars to $\widetilde{L}_m$ and identifying $\widetilde{L}_m\otimes_K H_i$ with $\widetilde{L}_m[N_i]$ via $\beta_i$, the scalar extension of the Hopf-Galois action of $H_i$ on $L_i$ coincides with the  action of $\widetilde{L}_m[N_i]$ on $\widetilde{L}_m\otimes_K L_i$.

\begin{lem}\label{lem: equivaction}
 Suppose that  for $n\leq m$ we have a surjective group homomorphism $q_{n,m}:N_m\rightarrow N_n$ such that its $\widetilde{L}_m$-extension $f_{n,m}:\widetilde{L}_m\otimes K[N_m]\longrightarrow \widetilde{L}_m\otimes K[N_n]$ descends to a morphism $\pi_{n,m}: H_m\rightarrow H_n$,  that is, we have $(1\otimes \widetilde{\beta}_n)\circ (1\otimes\pi_{n,m})=f_{n,m}\circ \beta_m$. Then the following are equivalent. 
 \begin{enumerate}
     \item\label{lemcond:1} For every $x\in X_m$ and $\eta\in N_m$ $\psi_{n,m}(x\cdot\eta)=\psi_{n,m}(x)\cdot q_{n,m}(\eta)$.
     \item\label{lemcond:2} For every $h\in H_m$ and $y\in L_n$, $h\cdot y=\pi_{n,m}(h)\cdot y$, where the equality is taken in $L_m$ via the inclusion $L_n\subseteq L_m$.
 \end{enumerate}
\end{lem}

\begin{proof}
 Let $ \iota_{n,m}: \widetilde{L}_m\otimes L_n \longrightarrow \widetilde{L}_m\otimes L_m $ be the map induced by the inclusion $L_n\subseteq L_m$. We first establish  the compatibility of $\phi_n$ and $\phi_m$ with the restriction map  $ \psi_{n,m}:X_m\longrightarrow X_n.$ For $a\in\widetilde{L}_m$, $y\in L_n$, and $x\in X_m$, we have  $$\begin{aligned} \phi_m\bigl(\iota_{n,m}(a\otimes y)\bigr)(x) &= \phi_m(a\otimes y)(x)\\
 &= a\,x(y)\\
 &= a\,\psi_{n,m}(x)(y)\quad\text{since $y\in L_n$}\\ 
 &= \phi_n(a\otimes y)\bigl(\psi_{n,m}(x)\bigr).
 \end{aligned} $$
Since $L_n\subseteq L_m$, the Hopf-Galois action $\mu_m:H_m\otimes L_m\rightarrow L_m$ restricts to $\mu_m:H_m\otimes L_n\rightarrow L_m$. By extending scalars, we have  $1\otimes\mu_m:(\widetilde{L}_m\otimes H_m)\otimes_{\widetilde{L}_m} (\widetilde{L}_m\otimes L_n)\rightarrow \widetilde{L}_m\otimes L_m$. Hence, we obtain the following diagram, where the horizontal rows of the right square diagram are the linear actions.
$$ \xymatrixcolsep{2.5pc}\xymatrix{
   (\widetilde{L}_m\otimes H_m)\otimes_{\widetilde{L}_m} (\widetilde{L}_m\otimes L_n) \ar[r]^-{\beta_m\otimes 1} \ar[d]_{(1\otimes\pi_{n,m})\otimes 1} & \widetilde{L}_m[N_m]\otimes_{\widetilde{L}_m}  \widetilde{L}_m\otimes L_n \ar[d]^{f_{n,m}\otimes 1}  \ar[r]^{ } &  \widetilde{L}_m\otimes L_m  \\          (\widetilde{L}_m\otimes H_n)\otimes_{\widetilde{L}_m} (\widetilde{L}_m\otimes L_n) \ar[r]_-{ \widetilde{\beta}_n\otimes 1 } &\widetilde{L}_m[N_n]\otimes_{\widetilde{L}_m}  \widetilde{L}_m\otimes L_n   \ar@{->>}[r]^{ }  & \ar@{^{(}->}[u]_{i_{n,m}}  \widetilde{L}_m\otimes L_n .}$$
 By hypothesis, the left square diagram  commutes. The first statement is equivalent to the commutativity of the right square diagram. Indeed, the right square commutes if and only if  for any $y\in L_n$ and $\eta\in N_m$,  \begin{align*}
    i_{n,m}(q_{n,m}(\eta)\cdot(1\otimes y))=\eta\cdot (1\otimes y)&   \Longleftrightarrow \forall x\in X_m,\; \phi_m\bigl(i_{n,m}(q_{n,m}(\eta)\cdot(1\otimes y))\bigr)(x)= \phi_m\bigl(\eta\cdot (1\otimes y)\bigr)(x)\\
    &\Longleftrightarrow \forall x\in X_m,\; \phi_n(q_{n,m}(\eta)\cdot(1\otimes y))(\psi_{n,m}(x))= (x\cdot\eta)(y)\\
    &\Longleftrightarrow \forall x\in X_m,\; (\psi_{n,m}(x)\cdot  q_{n,m}(\eta))(y)=\psi_{n,m}(x\cdot\eta)(y).
\end{align*}
 The second statement means that $\mu_m(h\otimes y)= \mu_n(\pi_{n,m}(h)\otimes y)$. After linear extension over $\widetilde{L}_m$ this is equivalent to $1\otimes \mu_m(h\otimes y)=i_{n,m}\big(1\otimes\mu_n(\pi_{n,m}(h)\otimes y) \big)$. Therefore, the second statement is equivalent to the commutativity of the large rectangle. This shows that the two conditions are equivalent.
\end{proof}

\begin{lem}\label{lem:fixed-finite-level}   Suppose that the two equivalent assertions of Lemma~\ref{lem: equivaction} hold.  Then $L_n$ is an $H_m$-normal sub-extension of $L_m$. More precisely, $$J_{H_m}(L_n):= \{h\in H_m\mid h\cdot y=0\text{ for every }y\in L_n\} = \ker(\pi_{n,m}).$$  Consequently,  $L_m^{\ker(\pi_{n,m})}=L_n.$ \end{lem} 

\begin{proof}  Let $h\in H_m$ and $y\in L_n$.  Since $\pi_{n,m}(h)\in H_n$, then  $\pi_{n,m}(h)\cdot y\in L_n.$ Hence $h\cdot y\in L_n$, and therefore $ H_m\cdot L_n\subseteq L_n.$ It is clear that $\ker(\pi_{n,m})\subseteq J_{H_m}(L_n)$.   Conversely,  let $h\in J_{H_m}(L_n)$. Then $h\cdot y=\pi_{n,m}(h)\cdot y=0$ for every $y\in L_n$. Hence $\pi_{n,m}(h)$  annihilates $L_n$ for every $h\in H_m$. As $L_n/K$ is $H_n$-Galois, the canonical map $\operatorname{\underline{can}}_n: L_n\otimes_KH_n\longrightarrow\operatorname{End}_K(L_n)$ is bijective. It follows that $\operatorname{\underline{can}}_n(1\otimes \pi_{n,m}(h))(y)=\pi_{n,m}(h)\cdot y=0$ implies that $1\otimes \pi_{n,m}(h)=0 $.  Thus $ \pi_{n,m}(h)=0$ and $h\in\ker\pi_{n,m}$.  Since $\pi_{n,m}$ is surjective we have $H_m/J_{H_m}(L_n)\cong H_n$. Thus $L_n$ is an $H_m$-normal sub-extension. Theorem~\ref{thm: corespo} gives $ L_m^{\ker(\pi_{n,m})}=L_m^{J_{H_m}(L_n)}=L_n$.  \end{proof}

\begin{proof}[Proof of Theorem \ref{thm: main}]
Let $n,m\in I$ be such that $n\leq m$. The   morphism $q_{n,m}$ induces a  morphism of $\widetilde{L}_m$-Hopf algebras $f_{n,m}:\widetilde{L}_m\otimes K[N_m]\longrightarrow \widetilde{L}_m\otimes K[N_n]$. Condition~(\ref{cond: 2}) and Proposition~\ref{pro: desct} tell us that $f_{n,m}$ descends to a unique  morphism of $K$-Hopf algebras  $\pi_{n,m}:H_m\longrightarrow H_{n}$. That is, the following diagram commutes.
$$\xymatrix{ \widetilde{L}_m\otimes H_m\ar[r]^{ 1\otimes \pi_{n,m}} \ar[d]_{ \beta_m} &   \widetilde{L}_m\otimes \widetilde{H}_n  \ar[d]^{\widetilde{\beta}_n =1\otimes \beta_n}   \\ \widetilde{L}_m[N_m] \ar[r]_-{ f_{n,m}} & \widetilde{L}_m[N_n].}$$
Moreover, since $f_{n,m}$ is surjective then $\pi_{n,m}$ is also surjective.
 Let $l,n,m\in I$ be such that  $l\leq n\leq m$.   Since $q_{l,m}=q_{l,n}\circ q_{n,m}$, we have $f_{l,m}=f_{l,n}\circ f_{n,m}$. By the  uniqueness of the descended morphism, we obtain $\pi_{l,m}=\pi_{l,n}\circ \pi_{n,m}$. Similarly, we have $\pi_{n,n}=\id_{H_n}$. Let $ H:=\varprojlim_n H_n$,  $p_n:H\rightarrow H_n,\;(h_m)_m\mapsto h_n$ be the natural projection and $J_n=\ker p_n$. Since all the $\pi_{n,m}$ are surjective for $n,m\in I$ with $n\leq m$, then $p_n$ is surjective for all $n\in I$. We have $p_n=\pi_{n,m}\circ p_m$ and $H/J_n\cong H_n$. Then the family of proartinian ideals $(J_n)_{n\in I}$ forms a basis of open neighbourhoods of $0$ so that $H$ is a proartinian Hopf algebra. We let $H$ act on $L$ as follows. For $h=(h_n)_n\in H$ and $y\in L$, we choose $n\in I$ such that $y\in L_n$ and define $h\cdot y:=p_n(h)\cdot y$. Condition~(\ref{cond: 3}) tells us that this action is well defined, that is, it does not depend on the choice of $n$. Indeed, suppose $y\in L_m$ for some $m\in I$. As $I$ is directed, there exists $l\in I$ such that $n\leq l$ and $m\leq l$. By Lemma~\ref{lem: equivaction}, $p_l(h)\cdot y = \pi_{n,l}(p_l(h))\cdot y = p_n(h)\cdot y$.  Similarly, $ p_l(h)\cdot y=p_m(h)\cdot y$.  Hence $ p_n(h)\cdot y=p_m(h)\cdot y$,  so the action is well defined.  
 Since for all  $n\in I$, $L_n$ is an $H_n$-module algebra, one checks that $L$ is a discrete $H$-module algebra. Now, we need to show that the canonical continuous  map $$\operatorname{\textbf{\underline{can}}}: L\widehat{\otimes} H\rightarrow \en_K(L),\; x\widehat{\otimes} h\mapsto (y\mapsto x(h\cdot y)) $$ is bijective.   Since  for any  $n\in I$, the finite extension  $L_n/K$ is Hopf Galois, the $K$-linear map $\operatorname{\underline{can}}_n: L_n \otimes H_n\rightarrow \en_K(L_n)$ is bijective.  After a scalar extension from  $L_n$ to $L$, as the map $L\otimes_{L_n}\en_K(L_n)\rightarrow \hm_K(L_n,L),\; x\otimes f\mapsto (y\mapsto xf(y))$ is an isomorphism,  we obtain a bijective morphism  $$L\otimes  \operatorname{\underline{can}}_n:L\otimes H_n\rightarrow \hm_K(L_n,L).$$ For any $n\leq m$, let $\operatorname{res}_{n,m}:\hm_K(L_m,L)\rightarrow \hm_K(L_n,L),\;f\mapsto f_{|L_n}$.  By the compatibility of the Hopf Galois actions, the following diagram commutes.
 $$ \xymatrixcolsep{3pc}\xymatrix{ L\otimes H_m\ar[r]^-{ L\otimes  \operatorname{\underline{can}}_m} \ar[d]_{1\otimes \pi_{n,m}} &   \hm_K(L_m,L)  \ar[d]^{\operatorname{res}_{n,m}}   \\ L\otimes H_n \ar[r]_-{ L\otimes\operatorname{\underline{can}}_n} & \hm_K(L_n,L).}$$
 Therefore, the family $(L\otimes  \operatorname{\underline{can}}_n)_{n\in I}$ defines  a morphism of inverse limits: $$\varprojlim_n (L\otimes  \operatorname{\underline{can}}_n):\varprojlim_n (L\otimes H_n)\longrightarrow \varprojlim_n \hm_K(L_n,L) $$ Since $\varprojlim_n (L\otimes H_n)=L\widehat{\otimes} H$ and $\varprojlim_n \hm_K(L_n,L)=\hm_K(\bigcup_n L_n, L)=\hm_K(L,L)=\en_K(L)$, then $$ \varprojlim_n (L\otimes  \operatorname{\underline{can}}_n):L\widehat{\otimes} H\longrightarrow \en_K(L).$$
 Moreover, since each $L\otimes  \operatorname{\underline{can}}_n$ is bijective, then  $\varprojlim_n (L\otimes  \operatorname{\underline{can}}_n)$ is also bijective and  by definition, it is the unique morphism that makes the following diagram commutative for all $n\in I$.$$\xymatrixcolsep{3.5pc}\xymatrix{
	L\widehat{\otimes} H  \ar[r]^-{ } \ar[d]_{1\otimes p_n} & \en_K(L) \ar[d]^{res_n}\\
	L\otimes H_n \ar[r]_-{\operatorname{L\otimes \underline{can}}_n } &\hm_K(L_n,L)}$$
To conclude the proof we need to show that $L^{J_n}=L_n$ so that 
  $\operatorname{\textbf{\underline{can}}}$ also makes the above diagram commutative, and thus $\operatorname{\textbf{\underline{can}}}=\varprojlim_n L\otimes  \operatorname{\underline{can}}_n$.
Recall that $ L^{J_n}:=\{y\in L\mid h\cdot y=0\; \forall h\in J_n\}$.   Let $y\in L_n$ and $h\in J_n$, we have  $h\cdot y=p_n(h)\cdot y=0$.  Hence $L_n \subseteq L^{J_n}$. Conversely, let $y\in L^{J_n}$. We choose $m\in I$ such that $y\in L_m$. Since $I$ is directed, after replacing $m$ by a common upper bound, we may assume that $n\leq m$. Since $p_m:H\to H_m$  is surjective and $p_n=\pi_{n,m}\circ p_m$, then  $p_m(J_n)=\ker(\pi_{n,m})$.   As $J_n$ annihilates $y$, it follows that $\ker(\pi_{n,m})$ annihilates $y$, because $y\in L_m$ and $J_n\cdot y=p_m(J_n)\cdot y$. Thus by Lemma~\ref{lem:fixed-finite-level} , $ y\in L_m^{\ker(\pi_{n,m})}=L_n$.   
\end{proof}

\begin{rem}
let $\widetilde{L}:=\bigcup_{n\in I}\widetilde{L}_n$, $G:=\gal(\widetilde{L}/K)=\varprojlim_n G_n$, $G':=\gal(\widetilde{L}/L)$, and $N:=\varprojlim_n N_n$. We have $$X:=G/G'\simeq\hm_{K\operatorname{-Alg}}(L,\widetilde{L})=\varprojlim_n\hm_{K\operatorname{-Alg}}(L_n,\widetilde{L})= \varprojlim_n\hm_{K\operatorname{-Alg}}(L_n,\widetilde{L}_n)=\varprojlim_nX_n.$$   Since for all $n\in I$, $G_n$ acts on $X_n$, we have a continuous action of $G$ on $X$. Moreover, since for all $n\in I$, $N_n$ is a regular subgroup of $\operatorname{Perm}(X_n)$ normalised by $G_n$, it follows that $N$ is a regular subgroup of $\operatorname{Perm}(X)$ normalised by $G$. 
\end{rem}

\section{Hopf Galois structure on $K(Q_n)/K$ }\label{section 3}
 In this section, $K$ denotes a number field. We first establish sufficient conditions for the extensions $K(Q_n)/K$ to admit Hopf-Galois structures. We then answer Question~\ref{question} negatively in general by constructing examples of such extensions that are not Hopf-Galois.

\subsection{Existence of a Hopf-Galois structure}\label{subsection 3.2}

We start this section by recalling the settings.  Let $E/K$ be an elliptic curve of rank at least $1$, and let $P\in E(K)$ be a fixed rational point of infinite order. For any positive integer $n$, let $Q_n\in E(\overline{K})$ be a point such that $nQ_n=P$. We consider the field $L_n:=K(Q_n),$ obtained by adjoining to $K$ the coordinates of $Q_n$. We denote by $\widetilde{L}_n$ its Galois closure over $K$ and $G_n=\gal(\widetilde{L}_n/K)$.   
The map  $$\delta_n:G_n\rightarrow E[n]\,, \sigma\mapsto  \sigma(Q_n)-Q_n$$ is a $1$-cocycle of $G_n$ with values in $U_n:=E(\widetilde{L}_n)[n]$. Although it is not a group homomorphism, we can still define its kernel as the subset of elements in $G_n$ whose image under $\delta_n $ is zero. It turns out that this kernel is a subgroup of $G_n$ which coincides with the stabiliser of the action of $G_n$ on $Q_n$. The latter is exactly the Galois group $G_n':=\gal(\widetilde{L}_n/L_n)$. Hence,  $G_n/G_n'$ is in bijection with $\im \delta_n$ and the orbit  $\operatorname{Orb}_{G_n}(Q_n)$  of $Q_n$ under $G_n$ is in bijection with the set $Q_n+\im\delta_n$. Thus, $X_n:=\hm_{K\operatorname{-Alg}}(L_n,\widetilde{L}_n)\simeq Q_n+\im\delta_n$.  Since the Galois conjugates of $Q_n$ are of the form $Q_n+\delta_{n}(\tau)$ for some $\tau\in G_n$, it follows that $\widetilde{L}_n=K(\im\delta_n, Q_n)$.  However, $\im\delta_n$ is not $G_n$-stable in general. We then consider the $G_n$-stable set  $$I_n=\{\sigma\delta_n(\tau)\;|\, \sigma,\tau\in G_n\}= \{\sigma\tau(Q_n)-\sigma(Q_n)\;| \sigma,\tau\in G_n\}$$ and  the field $K_n:=K(I_n)$ which is a Galois extension of $K$ contained in $\widetilde{L}_n$ and such that $K_nL_n=\widetilde{L}_n$.

 Since $I_n\subset U_n$,  the restriction of $\delta_n$ to the subgroup  $\Gamma_n:=\gal(\widetilde{L}_n/K_n)$ is an injective group homomorphism. Hence, the group $\Gamma_n$ is identified via $\delta_n$ with a subgroup of $U_n$, which we denote $V_n$.

The group $G_n$ acts naturally by field automorphisms on $U_n$. We obtain a group homomorphism $G_n\rightarrow \aut(U_n),\, \sigma\mapsto M_\sigma$. This induces an affine action of $G_n$ on $U_n$ given by $$\sigma\cdot u:= \delta_n(\sigma)+M_\sigma u,$$ for all $u\in U_n$ and $\sigma \in G_n$. Note that $\im\delta_n$ is stable under the affine action. Indeed, for $u=\delta_n(\tau)\in\im\delta_n$ with $\tau\in G_n$ and $\sigma\in G_n$, we have $\delta_n(\sigma)+M_\sigma u=\delta_n(\sigma\tau)\in\im\delta_n$. Thus the action of $G_n$ on $X_n$ by field automorphisms is given by:
$$ \sigma(Q_n+u)=Q_n+\delta_n(\sigma)+M_\sigma u=Q_n+\sigma\cdot u,$$ 
for all  $\sigma\in G_n,$ and $Q_n+u\in X_n $.

As $G_n/G_n'$ is in bijection with $X_n$, the injective  morphism $\lambda$ in Theorem~\ref{thm:G-P} can be viewed as $$ \lambda: G_n\longrightarrow \operatorname{Perm}(X_n),\,\sigma\mapsto\lambda(\sigma): Q_n+u\mapsto \sigma(Q_n+u)$$ and therefore the extension $L_n/K$ is Hopf-Galois if and only if there is a regular subgroup  of $\operatorname{Perm}(X_n)$ normalised by $G_n$. 

\begin{lem}\label{lem:unioncoset}
Let $G$ be a group and let $G_0$ be a subgroup of $G$. Let $A$ be an abelian group, and let $f:G\longrightarrow A$ be a map such that $f(hg)=f(h)+f(g)$ for all $h\in G_0$ and $g\in G$. Then $f(G)$ is a disjoint union of cosets of $f(G_0)$ in $A$; namely, $f(G)=\bigsqcup_{r\in \mathcal{R}}(f(G_0)+f(r))$, where $\mathcal{R}$ is a set of representatives for the left cosets of $G_0$ in $G$.
\end{lem}

\begin{proof}\label{lem:cosets}
We note that $f(G_0)$ is a subgroup of $A$. Every element $g\in G$ can be written as $g=hr$ for unique   $h\in G_0$ and  $r\in\mathcal R$. Hence $f(g)=f(h)+f(r)$, and so   $f(G)\subseteq \bigcup_{r\in\mathcal R}(f(G_0)+f(r))$. Conversely, if $a=f(h)+f(r)$ for some $h\in G_0$ and some $r\in\mathcal R$, by assumption, $a=f(hr)\in f(G)$.  Therefore,  $f(G)=\bigcup_{r\in\mathcal R}(f(G_0)+f(r))$. Finally, if two cosets $f(G_0)+a$ and $f(G_0)+b$ have non-empty intersection, then there exist $u,v\in f(G_0)$ such that $u+a=v+b$. Hence $a-b=v-u\in f(G_0)$, which implies $f(G_0)+a=f(G_0)+b$. Thus, the distinct cosets appearing in the above union are pairwise disjoint. Consequently, after choosing one representative for each distinct coset among the cosets $f(G_0)+f(r)$, $f(G)$ is a disjoint union of cosets of $f(G_0)$ in $A$.
\end{proof}

\begin{lem}\label{lem: keylem}
For any $\sigma\in G_n$ and $\gamma\in \Gamma_n$, we have $\delta_n(\sigma\gamma\sigma^{-1})=M_\sigma\delta_n(\gamma)$  and  $M_\sigma V_n= V_n$.
\end{lem}
\begin{proof}
Let $\sigma\in G_n$ and $\gamma\in \Gamma_n$. We have 
\begin{align*}
 \delta_n(\sigma\gamma\sigma^{-1})&=\delta_n(\sigma)+M_\sigma\delta_n(\gamma\sigma^{-1}) \\
 &=\delta_n(\sigma)+ M_\sigma(\delta_n(\gamma)+M_\gamma\delta_n(\sigma^{-1}))\\
 &=\delta_n(\sigma)+M_\sigma(\delta_n(\gamma)+\delta_n(\sigma^{-1}))\\
 &=\delta_n(\sigma)+M_\sigma \delta_n(\gamma) +M_\sigma\delta_n(\sigma^{-1})\\
 &= M_\sigma \delta_n(\gamma)\quad\text{ since }\quad \delta_n(\sigma^{-1})=-M_{\sigma^{-1}}\delta_n(\sigma)
\end{align*}
For the second part of the lemma, since all elements of $V_n$ are of the form $\delta_n(\gamma)$ for some $\gamma\in \Gamma_n$, and for a fixed $\sigma$, the map $\Gamma_n\longrightarrow \Gamma_n,\; \gamma\mapsto \sigma\gamma\sigma^{-1}$ is bijective, we have $M_\sigma V_n= V_n$.  
\end{proof}

Since $\delta_n(\gamma\sigma)=\delta_n(\gamma)+\delta_n(\sigma)$ for all $\gamma\in \Gamma_n$ and $\sigma\in G_n$, by Lemma~\ref{lem:unioncoset}, $\im\delta_n$ is a union of disjoint cosets of $V_n$. This implies that $ [L_n: K]=\#\operatorname{Orb}_{G_n}(Q_n)= \#\im\delta_n=r\cdot\#V_n$, where $r$ is the number of $V_n$-cosets in $\im\delta_n$. Note that the equivalence relation on $U_n$ that defines the quotient $U_n/V_n$ restricts to an equivalence relation on $\im\delta_n$.  We denote $\im\delta_n/V_n$ the set of these $r$ cosets.

\begin{pro}\label{prop: degree}
Let $F_n:=K_n\cap L_n$. Then $[F_n: K]=[G_n:\Gamma_nG_n']=r$. 
\end{pro}

\begin{proof}
 As  $\Gamma_n$ is normal in $G_n$, then $\Gamma_nG_n'$ is a subgroup of $G_n$ and the  map $ G_n\longrightarrow \im\delta_n/V_n,\, \sigma\mapsto \delta_n(\sigma)+V_n$  induces a well defined bijection  
$$\theta: G_n/\Gamma_nG_n'\longrightarrow \im\delta_n/V_n,\,\, \sigma \Gamma_nG_n'\mapsto \delta_n(\sigma)+V_n,$$ which gives the second equality. Indeed, if $\sigma_1\Gamma_nG_n'=\sigma_2\Gamma_nG_n'$, then we can write $\sigma_1=\sigma_2\gamma g$ where $\gamma\in \Gamma_n$ and $g\in G_n'$. We have 
\begin{align*}
\delta_n(\sigma_2\gamma g)&=\delta_n(\sigma_2)+M_{\sigma_2}(\delta_n(\gamma)+M_\gamma\delta_n(g))\\
&=\delta_n(\sigma_2)+M_{\sigma_2}\delta_n(\gamma)\quad \text{ since $\delta_n(g)=0$}\\
&= \delta_n(\sigma_2)+v \quad \text{ for some $v\in V_n$, since $M_{\sigma_2}V_n=V_n$.}
\end{align*}
It follows that $\delta_n(\sigma_1)+V_n=\delta_n(\sigma_2)+V_n$. Therefore, $\theta$ is well defined. The map is surjective by definition and for injectivity, if $\delta_n(\sigma_1)+V_n=\delta_n(\sigma_2)+V_n$, then $\delta_n(\sigma_1)-\delta_n(\sigma_2)\in V_n$ and 
\begin{align*}
 \delta_n(\sigma_2^{-1}\sigma_1)&=\delta_n(\sigma_2^{-1})+M_{\sigma_2^{-1}}\delta_n(\sigma_1)\\
 &=-M_{\sigma_2^{-1}}\delta_n(\sigma_2)+M_{\sigma_2^{-1}}\delta_n(\sigma_1)\\
 &=M_{\sigma_2^{-1}}(\delta_n(\sigma_1)-\delta_n(\sigma_2)).
\end{align*}
Hence, $\delta_n(\sigma_2^{-1}\sigma_1)\in V_n$. Let $\gamma\in \Gamma_n$ such that $\delta_n(\gamma)=\delta_n(\sigma_2^{-1}\sigma_1)$. With similar calculations as above, using the fact that $M_\gamma$ is trivial and that  $\delta_n$ is a group homomorphism on $\Gamma_n$, we get $\delta_n(\gamma^{-1}\sigma_2^{-1}\sigma_1)=0$. Thus  $\gamma^{-1}\sigma_2^{-1}\sigma_1\in G_n'$ which yields $\sigma_2^{-1}\sigma_1\in \Gamma_nG_n'$.  The first equality follows from the fact that $\widetilde{L}_n^{\Gamma_nG_n'}=\widetilde{L}_n^{\Gamma_n}\cap \widetilde{L}_n^{G_n'}=K_n\cap L_n=F_n$.
\end{proof}
\begin{pro}\label{prop: Im}
 If the image of $\delta_n$ is a subgroup of $U_n$, then the extension $L_n/K$ is Hopf-Galois
 for a Greither-Pareigis group isomorphic to $ \im\delta_n$.   
\end{pro}
 
\begin{proof}
Recall that $X_n\simeq Q_n+\im\delta_n$.  For each $u\in\im\delta_n$, we define $$t_u: X_n\longrightarrow X_n,\, Q_n+u'\mapsto
  Q_n+u'+u.$$ This is a well defined permutation of $X_n$ since $\im\delta_n$ is a subgroup of $U_n$.  The map $\im\delta_n\longrightarrow \operatorname{Perm}(X_n),\;u\mapsto t_u$ is an injective group homomorphism that induces a bijection between $\im\delta_n$ and its image $N:=\{t_u,\; u\in \im\delta_n\}\leq \operatorname{Perm}(X_n)$. Now, we will show that $N$ is a regular subgroup normalised  by $G_n$. Let $x=Q_n+u_1,\; y=Q_n+u_2\in X_n$. For $u=u_2-u_1$, we have $t_u(x)=y.$ Also, for any $u, u'\in \im\delta_n$, $t_u(Q_n+u')=Q_n+u'$  if and only if $u=0$. This shows that $N$ is a regular subgroup of $\operatorname{Perm}(X_n)$. Let $\sigma\in G_n$ and $t_u\in N$, we have for all $u'\in \im\delta_n$
\begin{align*}
\sigma t_u\sigma^{-1}(Q_n+u')&=\sigma t_u(Q_n+\delta_n(\sigma^{-1})+M_{\sigma^{-1}}u')\\
&= \sigma(Q_n+\delta_n(\sigma^{-1})+M_{\sigma^{-1}}u'+u)\\
&=Q_n+(u'+M_\sigma u)\\
&=t_{M_\sigma u}(Q_n+u')
\end{align*}
Thus, $\sigma t_u\sigma^{-1}=t_{M_\sigma u}\in N$. We conclude that $G_n$ normalises $N$.
\end{proof}
\begin{cor}\label{cor: m=1}
  If $K_n\cap L_n=K$, then $L_n/K$ is almost classically Galois  because $\im\delta_n\cong \Gamma_n$ is a subgroup of $U_n$.  Moreover, the corresponding  Greither-Pareigis group is isomorphic to $V_n$. 
\end{cor}

\begin{thm}\label{thm: quad}
 Suppose $[F_n:K]\leq 2$, then $L_n/K$ is Hopf-Galois for a Greither-Pareigis group isomophic to $V_n$ or $V_n\rtimes C_2$ .
\end{thm}

\begin{proof}
We know by Proposition~\ref{prop: degree} that $[F_n: K]=m$. The case $m=1$ is given by Corollary \ref{cor: m=1}. For  $m=2$, we can write   $\im\delta_n= V_n\cup (u_0+V_n)$, where $u_0\in \im\delta_n\setminus V_n$. Hence, $X_n=Q_n+V_n\cup (Q_n+u_0+V_n)$. For every $v\in V_n$, we define the map $$t_v: X_n\longrightarrow X_n,\, x\mapsto
  x+v.$$ Each $t_v$ is a well defined permutation of $X_n$. Indeed, let $x=Q_n+u$ with $u\in \im\delta_n$. Then $u=v'$ or $u=u_0+v'$ for some $v'\in V_n$. Therefore, $x+v\in Q_n+V_n$ or $x+h\in Q_n+u_0+V_n$. It is not difficult to see that $t_v$ is bijective.  We consider the group $T=\{t_v,\; v\in V_n\}\leq \operatorname{Perm}(X_n)$.   We define  the map $i:X_n\to X_n$ by
\[
i(Q_n+u)=Q_n+(u_0-u)
\]
This map is clearly an involution in $\operatorname{Perm}(X_n)$ since $u_0\notin V_n$. We define the group $$N:= \langle T, i\rangle \leq \operatorname{Perm}(X_n).$$
Since $i$ and $t_v$ do not commute for non trivial $t_v\in T$ (indeed $it_v=t_{-v}i$), and $T\cong V_n$, then $N\cong V_n\rtimes C_2 $.   We need to show that $N$ is a regular subgroup of $\operatorname{Perm}(X_n)$ normalised by $G_n$. Let $x=Q_n+v$ and $y=Q_n+u_0+v'$ for $v,v'\in V_n$, then $$t_{v+v'}i(y)=t_{v+v'}i(Q_n+u_0+v')=t_{v+v'}(Q_n-v')=Q_n+v=x.$$
The other cases for proving transitivity  are not difficult. Furthermore, the only element of $N$ fixing $Q_n$ is $t_0=\id$ since $t_vi(Q_n)\neq Q_n$ for all $v\in V_n$. Thus, $N$ is a regular subgroup of $\operatorname{Perm}(X_n)$. Using the same argument as in the proof of Proposition~\ref{prop: Im} one proves that $G_n$ normalises $T$. Next, we show that $G_n$ normalises $i$.

Let $\sigma\in G_n$, we have for all $u\in\im\delta_n$
  \begin{align*}
\sigma i\sigma^{-1}(Q_n+u)&= \sigma i(Q_n+\delta_n(\sigma^{-1})+M_{\sigma^{-1}}u)\\
&=\sigma(Q_n+(u_0-\delta_n(\sigma^{-1})-M_{\sigma^{-1}}u))\\
        &=Q_n+\delta_n(\sigma)+M_\sigma u_0-M_{\sigma}\delta_n(\sigma^{-1})-u\\
&=Q_n+ 2\delta_n(\sigma)+M_\sigma u_0-u.
\end{align*}
Let $v_\sigma:=M_\sigma u_0+ 2\delta_n(\sigma)-u_0\in U_n$, we have 
$$\sigma i\sigma^{-1}(Q_n+u)= t_{v_\sigma}i(Q_n+u).$$
Therefore, to conclude the proof, it remains to show that $v_\sigma\in V_n$. 
Let $\pi:U_n\longrightarrow U_n/V_n,\,u\mapsto \bar{u}.$ We have $\pi(\im\delta_n)=\{0,\bar{u_0}\}\subset U_n/V_n$.  Since $M_\sigma V_n=V_n$, $M_\sigma$ induces an automorphism $\overline{M_\sigma}\in \aut(U_n/V_n)$. Furthermore,  the affine map $$\varphi_\sigma: U_n\longrightarrow U_n, \, u\mapsto \delta_n(\sigma)+M_\sigma u $$ induces a well defined  map $$\overline{\varphi_\sigma}:U_n/V_n \longrightarrow U_n/V_n,\,\bar{u}\mapsto \overline{\delta_n(\sigma)}+\overline{M_\sigma}\bar{u}.$$  Since $\im\delta_n$ is stable under the affine action, $\varphi_\sigma(\im\delta_n)=\im\delta_n$. It follows that  $\overline{\varphi_\sigma}$ preserves $ \{0,\bar{u_0}\}$. Hence, $\overline{\delta_n(\sigma)}=\overline{\varphi_\sigma}(0)\in\{0,\bar{u_0}\}$ and $\overline{\delta_n(\sigma)}+\overline{M_\sigma}\bar{u_0}=\overline{\varphi_\sigma}(\bar{u_0})\in \{0,\bar{u_0}\} $.  Since $\overline{M_\sigma}\bar{u_0}\neq 0$,  the only possibilities are $$(\overline{\delta_n(\sigma)}, \overline{M_\sigma}\bar{u_0})=(0, \bar{u_0})\quad \text{and } \quad (\overline{\delta_n(\sigma)}, \overline{M_\sigma}\bar{u_0})=(\bar{u_0}, -\bar{u_0}). $$
Therefore, $\overline{v_\sigma}=\overline{M_\sigma}\bar{u_0}+2\overline{\delta_n(\sigma)}-\bar{u_0}=0$. Thus $v_\sigma\in V_n$.
\end{proof}

\begin{cor}
 For $n=3$, the extension $L_n/K$ is Hopf-Galois.
\end{cor} 
\begin{proof}
 Since $\#E[3]=9$, $\#H_3\in\{1,3,9\}$ and the only possible non trivial degree are $d\in \{3,6,9\}$.
 If $d=3,9$, then by the small degree result (see\cite[Theorem 4.6]{Greither}) and Proposition~\ref{prop: Im},   $L_3/K$ is Hopf-Galois. The case $d=6$ corresponds to $|H_3|=3$, and $m=2$ and the result follows  by Theorem~\ref{thm: quad}.  
\end{proof}

\subsection{Torsion representation and counter examples}\label{subsec: counter}

\begin{lem}\label{lem: counter}
Let $N$ be a group of order $72$ or $216$. Then $\#\aut(N)$ is not divisible by $5$. 
\end{lem}

\begin{proof}
We prove the result for $\#N=72$, the case $\#N=216$ is obtained with similar arguments. Suppose by contradiction that $5$ divides $\#\operatorname{Aut}(N)$. Then there exists an automorphism $\alpha\in \operatorname{Aut}(N)$ of order $5$. Let $\operatorname{Syl}_3(N)$ be the set of Sylow $3$-subgroups of $N$. Since $\#N=72=2^3\cdot 3^2,$ we have $\#\operatorname{Syl}_3(N)\in\{1,4\}$. The subgroup generated by $\alpha$ acts on $\operatorname{Syl}_3(N)$, and every orbit has size dividing $5$. Since $\#\operatorname{Syl}_3(N)<5$, every orbit has size $1$. Thus $\alpha$ fixes every Sylow $3$-subgroup of $N$. Let $H_3$ be a Sylow $3$-subgroup. Then $\alpha$ restricts to an automorphism of $H_3$ whose order divides $5$. But $H_3$ has order $9$, so $H_3$ is isomorphic either to $C_9$ or to $C_3\times C_3$. Their automorphism groups are respectively   $C_6$ and $\operatorname{GL}_2(\mathbb{F}_3)$. Neither has an element of order $5$. Hence $\alpha_{|H_3}=\operatorname{id}_{H_3}$. Similarly,   there exists a Sylow $2$-subgroup $H_2$ fixed by $\alpha$. Thus $\alpha$ restricts to an automorphism of $H_2$. Since $\#H_2=8$, the group $\operatorname{Aut}(H_2)$ has no element of order $5$. Therefore $\alpha_{|H_2}=\operatorname{id}_{H_2}$. Finally, since $H_2$ and $H_3$ have coprime orders and $\#H_2\#H_3=\#N$,
we have $N=H_2H_3$. Thus $\alpha$ is the identity on $N$, contradicting the fact that $\alpha$ has order $5$. 
\end{proof}

Let $E$, be an elliptic curve defined over $K$, and let $n$ be a positive integer. The group $E[n]$ of $n$-torsion points of $E$ is a free $\mathbb{Z}/n\mathbb{Z}$-module of rank $2$. After fixing a basis of $E[n]$, it is isomorphic to $(\mathbb{Z}/n\mathbb{Z})^2$ and $\aut(E[n])$ is isomorphic to $\gl_2(\mathbb{Z}/n\mathbb{Z})$. The absolute Galois group $G_K:=\gal(\overline{K}/K)$ acts $\mathbb{Z}/n\mathbb{Z}$-linearly on $E[n]$ and  this yields a group homomorphism $$\rho_{E,n}:G_K\longrightarrow \gl_2(\mathbb{Z}/n\mathbb{Z})$$
called the $\bmod$-$n$ torsion representation associated with the elliptic curve $E$. 

Now, suppose that $E$ has rank at least  $1$ and that $\rho_{E,\ell}$ is surjective for all primes $\ell$ outside a finite set of primes $D$. Let $p$ and $q$ be two distinct primes such that $p\notin D$, and set $n=pq$. Let $R\in E(K)$ be a point of infinite order, and let $S\in E(\overline{K})$ be such that $S\notin K_n$ and $qS=R$. Let $0\neq T\in E[p]$. We set $P=pR$ and $Q_n=S+T$. Then $nQ_n=pq(S+T)=pR=P$, so $Q_n$ is an $n$-division point of $P$. We claim that $L_n:=K(Q_n)=K(S,T)$. Indeed, since $Q_n=S+T$, we clearly have $K(Q_n)\subseteq K(S,T)$. Conversely, $qQ_n=qS+qT=R+qT$. Hence, $qT=qQ_n-R\in E(L_n)[p]$. Since multiplication by $q$ is invertible on $E[p]$, it follows that $T=(q^{-1}\bmod p)(qQ_n-R)\in E(L_n)$. Thus $S=Q_n-T\in E(L_n)$. Therefore, $K(S,T)\subseteq L_n$, and hence $L_n=K(S,T)$.

Let $\widetilde{L}_n$ be the Galois closure of $L_n/K$, and let $G_n:=\operatorname{Gal}(\widetilde{L}_n/K)$.

\begin{thm}\label{thm:degree}
With the notations as above, we have $$[L_{n}:K]\in\{(p^2-1)mq,\; 1\leq m\le q \}.$$ Moreover, $p$ divides the size of $G_n$.
\end{thm}

\begin{proof}
 By definition, $[K(T):K]=\#\operatorname{Orb}_{G_n}(T)$. Since  all the $G_n$-conjugates of $T$ are contained in $\widetilde{L}_n$ and  the mod-$p$ Galois representation is surjective then $\#\operatorname{Orb}_{G_n}(T)=p^2-1$. It follows that
$[L_n:K]=[L_n:K(T)][K(T):K]=(p^2-1)[L_n:K(T)]$. Let $M:=K(E[q],T)$. We denote by $J=\operatorname{Gal}(\widetilde{L}_n/M)$ and by  $G_T=\operatorname{Gal}(\widetilde{L}_n/K(T))$. We have $J\subseteq G_T$.  We consider the map
$$\delta^S:G_n\rightarrow E[q],\qquad \sigma\mapsto \delta_\sigma^S:=\sigma(S)-S.$$
It is not difficult to see that $\delta^S$ is a $1$-cocycle of $G_n$ with values in $E[q]$, and that its restriction to $J$ is an injective group homomorphism since the elements of $J$ are determined by their images at $S$. Therefore, $J$ is isomorphic to a subgroup $N_q$ of $E[q]$. We have $[L_n:K(T)]=\#\operatorname{Orb}_{G_T}(Q_n)=\#\operatorname{Orb}_{G_T}(S)=\#\{\delta_\sigma^S\mid \sigma\in G_T\}$.
The set $B:=\{\delta_\sigma^S\mid \sigma\in G_T\}$ is the image of the restriction of $\delta^S$ to $G_T$, and $\delta^S(\tau\sigma)=\delta^S(\tau)+\delta^S(\sigma)$ for all $\tau\in J$ and $\sigma\in G_T$. By Lemma~\ref{lem:cosets}, we can write $B$ as a union of disjoint cosets of $N_q$ in $E[q]$. If $m$ denotes the number of these cosets, then $\#B=m\cdot\#N_q$. Since $N_q$ is a non-trivial subgroup of $E[q]$, we have $\#N_q\in\{q,q^2\}$. Thus, $[L_n:K(T)]\in\{mq\mid 1\leq m\leq q\}$, and therefore, $[L_n:K]\in\{(p^2-1)mq,\;\mid 1\leq m\leq q
\}$. Moreover, as $G_n$ acts transitively on the non-zero elements of $E[p]$, it follows that all $p$-torsion points of $E$ are contained in $\widetilde{L}_n$. Hence $\operatorname{Gal}(K(E[p])/K)$ is a quotient of $G_n$. Since $\operatorname{Gal}(K(E[p])/K)\simeq \operatorname{GL}_2(\mathbb{Z}/p\mathbb{Z})$
and $p$ divides $\#\operatorname{GL}_2(\mathbb{Z}/p\mathbb{Z})$, we conclude that $p$ divides $\#G_n$.
\end{proof}

\begin{cor}\label{cor: counter}
Let $E$ be an elliptic curve of rank $1$ defined over a number field $K$ for which  $D=\{2,3\}$. Let $p=5$ and $q=3$. We have  $[L_n: K]\in \{72,144, 216\}$. In particular, if $[L_n: K]\in \{72, 216\}$, then $L_n/K$ is not Hopf-Galois.
 \end{cor}
\begin{proof}
If $L_n/K$ is Hopf-Galois, then Theorem~\ref{thm: Hol} tells us that $G_n\subseteq \operatorname{Hol}(N) $, where $N$ is a regular subgroup of the symmetric  group $S_{[L_n:K]}$. By Theorem~\ref{thm:degree},  $5$ divides $\#G_n$, but does not divide $\#\operatorname{Hol}(N)$ by Lemma~\ref{lem: counter}.
\end{proof}

An explicit example of Corollary \ref{cor: counter} is obtained with the  LMFDB curve 196a2 (See \cite{code} ).

\section{New examples of infinite Hopf-Galois extensions}\label{section 4}
In this section,  we prove Theorem~\ref{thm: main2}.  For the notation, we refer to Introduction and Section~\ref{subsection 3.2}. 
\begin{lem}\label{lem: inter}
Suppose that for every integer $n\geq 2$ we have $K_n\cap L_n=K$. Then for $n\mid m$, we have $\widetilde{L}_n \cap K_m= K_n$ and the restriction map $\Delta_{n,m}:\Gamma_m\longrightarrow \Gamma_n,\, \gamma\mapsto \gamma_{|\widetilde{L}_n}$ is a well defined surjective group homomorphism. 
\end{lem}
\begin{proof}
Since $K_m\cap L_m=K$, then $K_m\cap L_n=K$. It follows that $[K_mL_n: K_m]=[L_n:K]$. We also have that $[\widetilde{L}_n:K_n]=[L_n:K]$ and $K_m\widetilde{L}_n=K_mK_nL_n=K_mL_n$.  Hence, $[K_mL_n:K_m]=[\widetilde{L}_n: \widetilde{L}_n\cap K_m]$. Therefore, we obtain   $$[\widetilde{L}_n:K_n]=[L_n:K]=[K_mL_n:K_m]=[\widetilde{L}_n: \widetilde{L}_n\cap K_m].$$
Since $K_n\subseteq \widetilde{L}_n\cap K_m \subseteq \widetilde{L}_n$, we conclude that $\widetilde{L}_n \cap K_m= K_n$. The image of $\Delta_{n,m}$ is a  subgroup of $\gal(\widetilde{L}_n/\widetilde{L}_n\cap K_m)$. Its kernel is 
$\gal(\widetilde{L}_m/\widetilde{L}_nK_m)$ and so $\#\im\Delta_{n,m}=\frac{[\widetilde{L}_m:K_m]}{[\widetilde{L}_m: \widetilde{L}_nK_m]}=[\widetilde{L}_nK_m: K_m]=[K_mL_n:K_m]=[\widetilde{L}_n: \widetilde{L}_n\cap K_m]$.  Therefore, $\im\Delta_{n,m}=\gal(\widetilde{L}_n/\widetilde{L}_n\cap K_m)=\Gamma_n$
\end{proof}

\begin{thm}\label{thm: InfHopf}
Let $E/K$ be an elliptic curve of rank at least $1$ defined over $K$. Let $P \in E(K)$ be a point of infinite order. For any integers  $n,m\geq 2$, let $Q_n, Q_m\in E(\overline{K})$ such that $nQ_n = P$, and $\frac{m}{n}Q_m=Q_n$ if $n\mid m$.  Let $L_n := K(Q_n)$. Assume that  $K_n\cap L_n=K$ for all $n\ge 2$. Then the infinite extension $L_\infty := \bigcup_{n\geq 2} L_n$ of $K$ is almost classically Galois. 
\end{thm}

\begin{proof}
 Since $K_n\cap L_n=K$ for all $n\geq 2$, by Corollary~\ref{cor: m=1}, we know that the extension $L_n/K$ is almost classically Galois for a certain finite dimensional  $K$-Hopf algebra $H_n$ and the Greither-Pareigis group  is isomorphic to the  subgroup $V_n$ of $E[n]$. We need to show that the family of groups $(V_n)_{n\geq 2}$ satisfy Theorem~\ref{thm: main}.  For any integers $n$ and $m$ such that $n\mid m$, let $d=\frac{m}{n}$ and  let $q_{n,m}: V_m\rightarrow V_n,\;v\mapsto dv$.  By Lemma~\ref{lem: inter}, this is a well defined surjective group homomorphism. Indeed, write $v=\gamma(Q_m)-Q_m$ with $\gamma\in \Gamma_n$, then $dv=\gamma_{|\widetilde{L}_n}(Q_n)-Q_n\in V_n$. The surjectivity of $q_{n,m}$ follows from the surjectivity of $\Delta_{n,m}$. Let $\theta^n\in \operatorname{H}^1(G_n,\aut(V_n))$, $\theta^m\in \operatorname{H}^1(G_m,\aut(V_m))$ be the cocycles associated to $H_n$ and $H_m$ respectively. Note that for any $v=\gamma(Q_m)-Q_m\in V_m$, where $\gamma\in\Gamma_m$, using Lemma~\ref{lem: keylem}, we have $\theta^m_{\sigma}(v)=\sigma\gamma\sigma^{-1}(Q_m)-Q_m=\sigma(v)$. It is clear that $q_{n,n}=\id_{V_n}$ and $q_{l,m}=q_{l,n}\circ q_{n,m}$ whenever $l\leq n\leq m$. Let $v\in V_m,$   we have 
\begin{align*}
 q_{n,m}\theta^m_{\sigma}(v)&= d(\sigma(v))\\
 &= \sigma(dv)\\
 &=r_{n,m}(\sigma)(dv)\quad \text{ since $dv\in V_n$}\\
 &= \theta^n_{r_{n,m}(\sigma)}(dv)\\
 &=  \theta^n_{r_{n,m}(\sigma)}q_{n,m}(v).
\end{align*}

Now, we are left to show that for every $v\in V_m$ and every $x\in X_m$, we have  $ \psi_{n,m}(x\cdot v) =  \psi_{n,m}(x)\cdot q_{n,m}(v).$ Recall that $X_m:=\hm_{K\operatorname{-Alg}}(L_m,\widetilde{L}_m)\simeq Q_m+\im\delta_m=Q_m+V_m$. Under this identification, $V_m$ acts on $X_m$ by translation. More precisely, let $x\in X_m$ be any embedding. Its image is determined by the image of $Q_m$, which is of the form $Q_m+v'$ for some $v'\in V_m$. Hence,  for any $v\in V_m$, $(x\cdot v)(Q_m)= Q_m+v'+v$.  The restriction map  $\psi_{n,m}: X_m\rightarrow X_n $ is then induced by $q_{n,m}$, that is, $\psi_{n,m}(x)(Q_n)= Q_n+q_{n,m}(v')$, where $x: Q_m\mapsto Q_m+v'$. It follows that for any $v\in V_m$ and $x\in X_m$, we have
\begin{align*}
\psi_{n,m}(x\cdot v)(Q_n)&= Q_n+q_{n,m}(v'+v)\\
&=  \psi_{n,m}(x)(Q_n)+q_{n,m}(v)\\
&= (\psi_{n,m}(x)\cdot q_{n,m}(v))(Q_n)
\end{align*}
The extension $L_\infty/K$ is almost classically Galois because it is a direct limit of the finite classical Galois extensions $L_n/K$.
\end{proof}

As we saw in Section \ref{section 2} and \ref{section 3}, the extension $L_n/K$ is not always Hopf-Galois; the obstruction lies in $F_n=K_n\cap L_n$ which is not always equal to $K$.  

\begin{lem}
There exists a constant $C=C(E,P,K)$ depending only on $E$, $P$ and $K$ such that $[F_n:K]\leq C$ for all $n\geq 2$.
\end{lem}
\begin{proof}
By \cite[Theorem 1]{Bertrand}) there exists a constant $C=C(E,P,K)$ such that $\#U_n/V_n\leq C$. Using Proposition \ref{prop: degree}, we have $[F_n:K]=\#\im\delta_n/V_n\leq \#U_n/V_n\leq C $.
\end{proof}

Let $n,m $ be two integers such that $n\mid m$. Since  $L_n\subseteq L_m$ and  $K_n\subseteq K_m$ then $F_n\subseteq F_m$. Let  $K':=\bigcup_{n\geq 2}F_n$ and  $K_{\infty}:=\bigcup_{n\geq 2}K_n$.

\begin{pro}\label{pro: UniHopf}
The field $K'=\bigcup_{n\geq 2}F_n$ is a finite extension of  $K$ and for any  $n\geq 2$, $K'=K' K_n\cap K' L_n$.
\end{pro}
\begin{proof}
Let us consider any $C+1$ elements $x_1,x_n,\ldots x_{C+1}\in K'$. Then for every $i\in\{1,\ldots C+1\}$, there exists an integer $n_i\geq 2$ such that $x_{i}\in F_{n_i}$. It follows that $x_1,x_n,\ldots x_{C+1}\in F_{n_0}$, where $n_0=\operatorname{lcm}(n_1,\ldots n_{C+1})$. Since $[F_{n_0}:K]\leq C$, we deduce that the elements $x_1,x_n,\ldots x_{C+1}$ are linearly dependent over $K$. Therefore, any family of $C+1$ elements of $K'$ is linearly dependent over $K$. Thus $[K':K]\leq C$.  Moreover, it is straightforward that $K'=K_\infty\cap L_\infty$ and then for all $n\in\mathbb{N}$, $K' K_n\cap K' L_n\subseteq K_\infty\cap L_\infty=K'$. The fact that $K'\subset  K' K_n\cap K' L_n$ gives us the equality.
\end{proof}

\begin{proof}[Proof of Theorem \ref{thm: main2}]
By Proposition~\ref{pro: UniHopf} we have $K' K_n\cap K' L_n=K'$. Hence, by Corollary~\ref{cor: m=1} the extension $K'L_n/K'$ is Hopf-Galois for all $n\geq 2$. After replacing the base field $K$ by $K'$, and for every integer $n\geq 2$, the fields $K_n$ $L_n$ and $\widetilde{L}_n$ are replaced by $K'K_n$, $K'L_n$ and $K'\widetilde{L}_n$ respectively. Applying  Theorem~\ref{thm: InfHopf} gives the result. 
\end{proof}

\begin{rem}
The techniques developed in Sections~\ref{section 3} and~\ref{section 4} can also be adapted to abelian varieties defined over number fields. They admit a similar adaptation to radical extensions $K(\sqrt[n]{a})$, where $a\in K^{\times}.$
\end{rem}

\section*{Acknowledgments}
I would like to thank my supervisor Gabor Wiese for suggesting this problem, for his guidance throughout its solution, and for his careful reading of all previous versions of the manuscript.


\let\cleardoublepage\clearpage

\bibliographystyle{plainurl}
\bibliography{mybib}

\bigskip

\noindent 
Leolin Nkuete\\
Department of Mathematics\\
University of Luxembourg \\
Luxembourg \\
\url{leolin.nkuete@uni.lu}

\end{document}